# A Center-Point Algorithm for Unit Commitment with Carbon Emission Trading

Linfeng Yang, *Member*, *IEEE*, Wei Li, Guo Chen, *Member*, *IEEE*, Beihua Fang, Chunming Tang, Zhaoyang Dong, *Fellow*, *IEEE*

*Abstract*—This paper proposes a global optimization method for it ensures finding good solutions while solving the unit commitment (UC) problem with carbon emission trading (CET). This method consists of two parts. In the first part, a sequence of linear integer-relaxed subproblems are first solved to rapidly generate a tight linear relaxation of the original mixed integer nonlinear programming problem (MINLP) model. In the second part, the algorithm introduces the idea of center-cut so that it can quickly find good solutions. The approach tested on 10 test instances with units ranging from 35 to 1560 over a scheduling period of 24h, and compared with state-of-the-art solver CPLEX. The results show that the proposed algorithm can find better solutions than CPLEX in a short time. And it is more suitable to solve large scale UC problem than CPLEX.

## NOMENCLATURE

Operator:
$\{\cdot\}^+$    max(0,·). When this operator is applied to a vector, it is performed element by element.

Indices:
$i$    Index for unit.
$t$    Index for time period.

Constants:
$N$    Total number of units.
$T$    Total number of time periods.
$\alpha_i, \beta_i, \gamma_i$    Coefficients of the quadratic production cost function of unit $i$.
$C_{\text{cold},i}$    Cold startup cost of unit $i$.
$C_{\text{hot},i}$    Hot startup cost of unit $i$.
$T_{\text{cold},i}$    Cold startup time of unit $i$.
$\underline{P}_i$    Minimum power output of unit $i$.
$\overline{P}_i$    Maximum power output of unit $i$.
$P_{\text{D},t}$    System load demand in period $t$.
$R_t$    Spinning reserve requirement in period $t$.
$P_{\text{up},i}$    Ramp up limit of unit $i$.
$P_{\text{down},i}$    Ramp down limit of unit $i$.
$P_{\text{start},i}$    Startup ramp limit of unit $i$.
$P_{\text{shut},i}$    Shutdown ramp limit of unit $i$.
$u_{i,0}$    Initial commitment state of unit $i$.
$T_{i,0}$    Number of periods for which unit $i$ has been online ($+$) or offline ($-$) prior to the first period of the time span (end of period 0).
$\underline{T}_{\text{on},i}$    Minimum up time of unit $i$.
$\underline{T}_{\text{off},i}$    Minimum down time of unit $i$.
$U_i$    $\{\min[T, u_{i,0}(\underline{T}_{\text{on},i} - T_{i,0})]\}^+$
$L_i$    $\{\min[T, (1 - u_{i,0})(\underline{T}_{\text{off},i} + T_{i,0})]\}^+$.
$\pi_{\text{b}}$    Price of emission allowances bought on the market.
$\pi_{\text{s}}$    Price of emission allowances sold on the market.
$E_0$    Total emission allowances of $CO_2$.

Variables:
$u_{i,t}$    On/off status of unit $i$ in period $t$.
$s_{i,t}$    Startup status of unit $i$ in period $t$.
$P_{i,t}$    Power output of unit $i$ in period $t$.
$S_{i,t}$    Startup cost of unit $i$ in period $t$.
$u_{\text{b}}$    Binary variable that is equal to 1 if the emission allowances bought on the market and 0 means no buying.
$u_{\text{s}}$    It is 1 if the emission allowances sold on the market and 0 means no selling.
$\Delta E_{\text{b}}$    Emission allowances bought on the market.
$\Delta E_{\text{s}}$    Emission allowances sold on the market.

## I. INTRODUCTION

In the power system, it is crucial to develop a production planning of generating units based on actual conditions and requirements. Therefore, a class of mathematical models have been developed, which are referred to the Unit Commitment (UC) problem. Generally, the purpose of UC problem is to minimize energy production costs while meeting the electricity demand and reserve requirements and different operating constraints of the generating units. In recent years, as the public pays more and more attention to the greenhouse effect, some researchers have begun to study UC problem considering carbon emissions [1]-[2].

As a basic problem in the power system, UC problem always can be formulated as a mixed-integer nonlinear programming (MINLP) problem, which is very challenging to solve efficiently due to the complexity of the problem. For last several decades, the UC problem has been the focus of significant research. There are many (meta-)heuristic and deterministic

This work was supported by the Natural Science Foundation of China (51767003, 71861002, 61862004), the Guangxi Natural Science Foundation (2017GXNSFBA198238). And This research is partially supported by Australian Research Council under grants DP180103217 and IH180100020.

L.F. Yang is with the School of Computer Electronics and Information, Guangxi University, Nanning 530004, China, and also with the Guangxi Key Laboratory of Multimedia Communication and Network Technology, Guangxi University, China. (e-mail: ylf@gxu.edu.cn)

W. Li and B.H. Fang are with the School of Computer Electronics and Information, Guangxi University, Nanning 530004, China. (e-mail: 2863053347@qq.com; 747751326@qq.com)

G. Chen and Z.Y. Dong are with the School of Electrical Engineering and Telecommunications, The University of NSW, Sydney, NSW 2052, Australia. (e-mail: andyguochen@gmail.com; zydong@ieee.org)

C.M. Tang is with the School of Mathematics and Information Science, Guangxi University, Nanning 530004, China. (e-mail: cmtang@gxu.edu.cn)

methods developed to solve it. (Meta-)Heuristic methods include priority list (PL) [3]-[5] and artificial intelligence (AI) algorithms [6]-[11]. However, the solutions obtained by (meta-)heuristic methods are difficult to guarantee good quality. Another class of techniques applied to UC problem are deterministic methods, which include dynamic programming (DP) [12]-[14], branch and bound (BB) [15], mixed-integer linear programming (MILP) [16]-[21], Lagrangian relaxation (LR) [22], [23]. For large scale problems, most deterministic methods require a long computation time to get high-quality solutions.

In this context, the main focus of this paper is the day-ahead operational planning of generating units considering carbon emission trading (CET). By setting a price for carbon dioxide emissions and limiting the carbon emission quota, the government can use economic means to protect the environment. In order to reduce economic costs, companies will actively reduce their emissions to operate beneath the cap, or they will spontaneously trade emissions rights with other companies. Therefore, CET is also called cap-and-trade. In general, the UC problem considering CET (UC-CET) can be formulated as a mixed-integer quadratic constrained programming (MIQCP) problem. For the convex MINLP problem, most deterministic methods are relied on the available well-developed theory, e.g. cutting plane (CP) method [24]. These methods include extended cutting plane (ECP) method [25], outer approximation method (OAM) [26]-[28], Benders decomposition (BD) [29], extended supporting hyperplane (ESH) algorithm [30], Center-Cut (CC) algorithm [31]. Among these methods, ECP is a variant of Kelley's cutting plane method. As one of the classical approaches for solving MINLP problem, the idea of ECP method has been used for reference by other linear approximation algorithms, such as OAM, BD, etc. But it tends to be unstable and converge slowly [25]. The OAM was originally proposed by Duran and Grossmann to solve convex MINLP [26]. Subsequently, some scholars proposed some variants of OAM and applied them to power system [27], [28]. BD and OAM have similar decomposition strategies. However, the difference between the two methods is that BD method is based on the dual information, while the OA method is based on the original information. Like OAM, BD method has also been applied to solve UC problems [29]. But OAM and BD usually take a long time to find feasible solutions because these two methods need to solve a MILP problem in each iteration [27], [29]. ESH algorithm looks for the linearization points on the boundary of the nonlinear constraint unlike ECP algorithm, where the linearization points usually violate nonlinear constraints [25], [30]. Meanwhile, unlike OAM and BD, where the linearization points are obtained by solving nonlinear programming (NLP) problems, ESH algorithm can get the linearization points by a line search procedure [30]. Therefore, ESH algorithm can quickly generate a linear approximation of the original problem. But ESH algorithm can only get a feasible solution after the end of the algorithm. Like ECP, OAM, BD and ESH, CC algorithm also needs to construct an approximate polyhedron of the original problem. However, the difference is that the algorithm chooses the Chebyshev Center of the polyhedron as the trail solution. It is this feature that makes CC algorithm only need a few iterations to find feasible solutions [31]. But CC algorithm usually takes a long time to construct a linear approximation of original MINLP problem. In [32], more details about the Chebyshev Center can be obtained.

In this paper, we propose a deterministic global optimization algorithm that can quickly find high-quality feasible solutions for UC-CET. The proposed algorithm can be divided into two steps. In the first step, a tight linear approximation of UC-CET can be rapidly generated by iteratively solving a sequence of linear continuous-relaxed subproblems, conducting linear search and constructing perspective-cuts. In the second step, high-quality solutions can be obtained by iteratively finding of the current linear approximation, searching the neighborhood of these solutions and adding new linear constraints. We discover that the CP algorithm can find better solutions by introducing the integer ellipse center.

The major contributions of this paper are summarized as follows:

(1) For UC-CET, a deterministic global optimization algorithm is proposed to find high-quality solutions faster than CPLEX.
(2) Compared to piecewise linear techniques, the CP algorithm can generate a tighter linear approximation with fewer number of cutting planes.
(3) Introducing the integer ellipse center innovatively, so that the performance of the CP algorithm is improved.

## II. MATHEMATICAL FORMULATION

In this section, we present the UC-CET model including commitment of generator units [18] and carbon emission trading [2].

### A. Objective Function

The primary purpose of production planning for UC-CET is to minimize the total cost of generating power,
$$F_{TC} = F_{TH} + F_{CET},$$
where $F_{TC}$ represents the total cost, $F_{TH}$ is the energy production cost of thermal units and $F_{CET}$ is the cost of carbon emissions trading.

*1) The energy production cost of thermal units ($F_{TH}$)*

Base on perspective-cut approximation [18], [33], the total cost of thermal units can be formulated as
$$F_{TH} = \sum_{i=1}^{N} \sum_{t=1}^{T} [z_{i,t} + C_{\text{hot},i} s_{i,t} + \tilde{S}_{i,t}] \quad (1)$$
where $z_{i,t}$ and $\tilde{S}_{i,t}$ are auxiliary variables constrained by the following constraints respectively,
$$z_{i,t} \geq \left(\frac{2\tilde{\gamma}_i l}{L} + \tilde{\beta}_i\right) \tilde{P}_{i,t} + \left(\tilde{\alpha}_i - \tilde{\gamma}_i \left(\frac{l}{L}\right)^2\right) u_{i,t}, \; l = 0,1,2,\dots,L \quad (2)$$
where $L$ is a given parameter, $\tilde{\alpha}_i = \alpha_i + \beta_i \underline{P}_i + \gamma_i (\underline{P}_i)^2$, $\tilde{\beta}_i = (\overline{P}_i - \underline{P}_i)(\beta_i + 2\gamma_i \underline{P}_i)$, $\tilde{\gamma}_i = \gamma_i (\overline{P}_i - \underline{P}_i)^2$, $\tilde{P}_{i,t} = \frac{P_{i,t} - u_{i,t} \underline{P}_i}{(\overline{P}_i - \underline{P}_i)}$.

$$\begin{cases} \tilde{S}_{i,t} \geq 0 \\ \tilde{S}_{i,t} \geq (C_{\text{cold},i} - C_{\text{hot},i}) \left[ s_{i,t} - \sum_{\tau = \max\{t - \underline{T}_{\text{off},i} - T_{\text{cold},i} - 1, 1\}}^{t-1} u_{i,\tau} - f_{\text{init},i,t} \right] \end{cases} \quad (3)$$

where $f_{\text{init},i,t} = 1$ when $t - \underline{T}_{\text{off},i} - T_{\text{cold},i} - 1 \leq 0$ and

$\{-T_{i,0}\}^+ < |t - \underline{T}_{\text{off},i} - T_{\text{cold},i} - 1| + 1$, $f_{\text{init},i,t} = 0$ otherwise.

*2) Cost of CET ($F_{CET}$)*

In this paper, the cost of CET can be formulated as
$$F_{\text{CET}} = \pi_b \Delta E_b - \pi_s \Delta E_s. \quad (4)$$

### B. Constraints

*1) Constraints of the thermal system*
- Unit generation limits:
$$0 \leq \tilde{P}_{i,t} \leq u_{i,t}. \quad (5)$$
- Power balance constraint:
$$\sum_{i=1}^{N}[\tilde{P}_{i,t}(\overline{P}_i - \underline{P}_i) + u_{i,t}\underline{P}_i] - P_{D,t} = 0. \quad (6)$$
- System spinning reserve requirement:
$$\sum_{i=1}^{N} u_{i,t}\overline{P}_i \geq P_{D,t} + R_t. \quad (7)$$
- Ramp rate limits:
$$\tilde{P}_{i,t} - \tilde{P}_{i,t-1} \leq u_{i,t}\tilde{P}_{\text{up},i} + s_{i,t}(\tilde{P}_{\text{start},i} - \tilde{P}_{\text{up},i}) \quad (8)$$
$$\tilde{P}_{i,t-1} - \tilde{P}_{i,t} \leq u_{i,t-1}\tilde{P}_{\text{shut},i} + (s_{i,t} - u_{i,t})(\tilde{P}_{\text{shut},i} - \tilde{P}_{\text{down},i}) \quad (9)$$

where $\tilde{P}_{\text{up},i} = \frac{P_{\text{up},i}}{\overline{P}_i - \underline{P}_i}$, $\tilde{P}_{\text{down},i} = \frac{P_{\text{down},i}}{\overline{P}_i - \underline{P}_i}$, $\tilde{P}_{\text{start},i} = \frac{P_{\text{start},i} - \underline{P}_i}{\overline{P}_i - \underline{P}_i}$, $\tilde{P}_{\text{shut},i} = \frac{P_{\text{shut},i} - \underline{P}_i}{\overline{P}_i - \underline{P}_i}$.

- Minimum up/down time constraints:
$$\sum_{\varpi=\{t-\underline{T}_{\text{on},i}\}^+ + 1}^{t} s_{i,\varpi} \leq u_{i,t}, \ t \in [U_i + 1, \dots, T] \quad (10)$$
$$\sum_{\varpi=\{t-\underline{T}_{\text{off},i}\}^+ + 1}^{t} s_{i,\varpi} \leq 1 - u_{i,\{t-\underline{T}_{\text{off},i}\}^+}, \ t \in [L_i + 1, \dots, T]. \quad (11)$$
- Initial status of units:
$$u_{i,t} = u_{i,0}, \ t \in [1, \dots, U_i + L_i]. \quad (12)$$

State constraints:
$$u_{i,t} - u_{i,t-1} \leq s_{i,t}. \quad (13)$$

*2) Constraints of CET*
$$\sum_{i=1}^{N}\sum_{t=1}^{T} u_{i,t} E_i(\tilde{P}_{i,t}) \leq E_0 + \Delta E_b - \Delta E_s \quad (14)$$
$$0 \leq \Delta E_b \leq u_b \Delta E_b^{\max} \quad (15)$$
$$0 \leq \Delta E_s \leq u_s \Delta E_s^{\max} \quad (16)$$
$$u_b + u_s \leq 1 \quad (17)$$

where $E_i(\tilde{P}_{i,t}) = \tilde{a}_i u_{i,t} + \tilde{b}_i \tilde{P}_{i,t} + \tilde{c}_i(\tilde{P}_{i,t})^2$, $\tilde{a}_i = a_i + b_i \underline{P}_i + c_i(\underline{P}_i)^2$, $\tilde{b}_i = (\overline{P}_i - \underline{P}_i)(b_i + 2c_i\underline{P}_i)$, $\tilde{c}_i = c_i(\overline{P}_i - \underline{P}_i)^2$. More details about projection substitution of (14) can be referred to the appendix A of reference [18]. And (14) restricts the total quantity of carbon emissions, (15) and (16) limit the maximum carbon emission quota can be bought or sold. (17) ensures that the selling and buying of carbon emission quota can be done asynchronously.

## III. CENTER-POINT ALGORITHM FOR UC-CET

In this section, we present our proposed algorithm for UC-CET. Before presenting the algorithm, we need to reformulate the UC-CET model, redefine the expression of the variables and introduce the difference between the perspective-cut and the general tangent.

### A. Model of the UC-CET

In this paper, the UC-CET can be formulated as a MIQCP problem.

$$\min F_{\text{TC}} = \sum_{i=1}^{N}\sum_{t=1}^{T}[z_{i,t} + C_{\text{hot},i}s_{i,t} + \tilde{S}_{i,t}] + \pi_b \Delta E_b - \pi_s \Delta E_s$$
$$s.t. \begin{cases} (2)(3)(5)-(17) \\ u_{i,t}, s_{i,t}, u_b, u_s \in \{0,1\}. \end{cases} \quad (18)$$

For constraint (14), introducing auxiliary variable $\eta_{\text{CET}}$, and let $\mathcal{G}(\chi) \coloneqq \sum_{i=1}^{N}\sum_{t=1}^{T} \tilde{c}_i(\tilde{P}_{i,t})^2 - \eta_{\text{CET}}$, we have
$$\sum_{i=1}^{N}\sum_{t=1}^{T}[\tilde{a}_i u_{i,t} + \tilde{b}_i \tilde{P}_{i,t}] + \eta_{\text{CET}} \leq E_0 + \Delta E_b - \Delta E_s \quad (19)$$
$$\mathcal{G}(\chi) \leq 0. \quad (20)$$

For the sake of convenience, we let $\chi = (u; P; x)$ where $u = (u_{i,t}; s_{i,t}; u_b; u_s)$, $P = (\tilde{P}_{i,t})$, $x = (\tilde{S}_{i,t}; \Delta E_b; \Delta E_s; \eta_{\text{CET}})$, $i = 1, \dots, N, t = 1, \dots, T$. Let $\ell(\chi) \coloneqq \sum_{i=1}^{N}\sum_{t=1}^{T}[z_{i,t} + C_{\text{hot},i}s_{i,t} + \tilde{S}_{i,t}] + \pi_b \Delta E_b - \pi_s \Delta E_s$. And denote $X_{\text{NL}} = \{\chi | (20)\}$, $X_{\text{L}} = \{\chi | (2)(3)(5)-(13)(15)-(17)(19)\}$, $X_{\text{I}} = \{\chi | u_{i,t}, s_{i,t}, u_b, u_s \in \{0,1\}\}$.

Finally the problem (18) can be equivalently rewritten as
$$\min_{\chi \in X_{\text{NL}} \cap X_{\text{L}} \cap X_{\text{I}}} \ell(\chi). \quad (21)$$

### B. Two linear approximations for $X_{NL}$

In this subsection, we will give two linear approximations for $X_{\text{NL}}$.

*1) Linear Approximation based on tangent*

Give a set of points $\Omega = \{\hat{\chi}^1, \hat{\chi}^2, \dots \hat{\chi}^\hbar\}$, we can generate a polyhedral outer approximation (OA) for $X_{\text{NL}}$ by constructing tangents at points of $\Omega$,
$$X_{\text{NL}}^{\text{TR}}(\Omega) = \left\{\chi \middle| \begin{array}{l} \sum_{i=1}^{N}\sum_{t=1}^{T}[2c_i\hat{\tilde{P}}_{i,t}^s\tilde{P}_{i,t} - c_i(\hat{\tilde{P}}_{i,t}^s)^2] - \eta_{\text{CET}} \leq 0 \\ s = 1, \dots, \hbar \end{array}\right\}. \quad (22)$$

*2) Linear Approximation based on perspective-cut*

With $0 \leq \tilde{P}_{i,t} \leq u_{i,t}$ and $u_{i,t} \in \{0,1\}$, we have the tight perspective-cuts linear approximation for $X_{\text{NL}}$,
$$X_{\text{NL}}^{\text{PR}}(\Omega) = \left\{\chi \middle| \begin{array}{l} \sum_{i=1}^{N}\sum_{t=1}^{T}[2c_i\hat{\tilde{P}}_{i,t}^s\tilde{P}_{i,t} - c_i u_{i,t}(\hat{\tilde{P}}_{i,t}^s)^2] - \eta_{\text{CET}} \leq 0 \\ s = 1, \dots, \hbar \end{array}\right\}. \quad (23)$$

**Theorem 1** Suppose $0 \leq \tilde{P}_{i,t} \leq u_{i,t}$ and $\eta_{\text{CET}} \geq 0$. Then $X_{\text{NL}}^{\text{PR}}(\Omega) = X_{\text{NL}}^{\text{TR}}(\Omega)$ for $u_{i,t} \in \{0,1\}$ and $X_{\text{NL}}^{\text{PR}}(\Omega) \subseteq X_{\text{NL}}^{\text{TR}}(\Omega) \supseteq X_{\text{NL}}$ for $u_{i,t} \in [0,1]$.

The equation can be verified by set $u_{i,t} = 0$ and $u_{i,t} = 1$ respectively.

For any $\hat{\tilde{P}}_{i,t}^s$, we have $2c_i\hat{\tilde{P}}_{i,t}^s\tilde{P}_{i,t} - c_i(\hat{\tilde{P}}_{i,t}^s)^2 \leq \tilde{c}_i(\tilde{P}_{i,t})^2$ because that $\tilde{c}_i(\tilde{P}_{i,t})^2$ is convex function. Then, the second inclusion relation can be obtained by summing. The first inclusion relation can be proved followed by $c_i u_{i,t}(\hat{\tilde{P}}_{i,t}^s)^2 \leq c_i(\hat{\tilde{P}}_{i,t}^s)^2$ [33].

Theorem 1 shows that $X_{\text{NL}}^{\text{PR}}(\Omega)$ is a tighter outer approximation set for $X_{\text{NL}}$ than $X_{\text{NL}}^{\text{TR}}(\Omega)$.

### C. Center-Point Method for UC-CET

The proposed center-point algorithm can be divided into two parts. The first part is mainly to linearly approximate $X_{\text{NL}}$ to obtain a tight approximate polyhedron, and the second part is mainly to find high-quality solutions.

*1) Generate a good linear approximation for $X_{NL}$*

• The center point step

In this step, we want to find an interior point of the nonlinear constraints by solving the NLP problem (24). Note that this point is not necessarily integer feasible.
$$\chi_{cp}^{NLP} = \arg \min_{\chi \in X_L} g(\chi). \tag{24}$$

It is obvious that $X_{NL}$ is convex and has interior points, then, suppose that the original UC-CET is feasible, NLP problem (24) has solutions and $g(\chi_{cp}^{NLP}) < 0$.

• The linear programming (LP) step

Now, we aim at quickly building and refining a polyhedral outer approximation of $X_{NL}$ based on linear approximation in each iteration. For giving points set $\Omega^R = \{\hat{\chi}^1, \hat{\chi}^2, \dots \hat{\chi}^\hbar\}$, we denote $X_{NL}^R(\Omega^R)$ as any polyhedral outer approximation of $X_{NL}$ corresponding to $\Omega^R$. We note that $X_{NL}^R(\Omega^R)$ can be $X_{NL}^{TR}(\Omega^R)$ or $X_{NL}^{PR}(\Omega^R)$ presented in section III.B.

After the solution of (24) is obtained, then the relaxation of (21) only considering the linear constraints is solved:
$$\chi^{LP} = \arg \min_{\chi \in X_{NL}^R(\Omega) \cap X_L} \ell(\chi). \tag{25}$$

If $g(\chi^{LP}) \leq 0$, $X_{NL}^R(\Omega^R)$ is a good linear approximation for $X_{NL}$. Otherwise, i.e. $g(\chi^{LP}) > 0$, a line search is performed between $\chi_{cp}^{NLP}$ and $\chi^{LP}$, i.e., the equation
$$\hat{\chi} = \lambda \chi_{cp}^{NLP} + (1 - \lambda) \chi^{LP} \tag{26}$$
can be used to find the value of $\lambda \in [0,1]$ such that $g(\hat{\chi}) = 0$. Then let $\Omega^R = \Omega^R \cup \{\hat{\chi}\}$ and update $X_{NL}^R(\Omega^R)$. Since the nonlinear constraint $g(\chi)$ is a quadratic function, the explicit expression of $\lambda$ can be directly given. The problem (25) is solved repeatedly until the iteration counter is greater than a maximum number, or until $g(\chi_{LP}) \leq 0$. By iteratively solving (25) and conducting line searches, $X_{NL}^R(\Omega^R)$ can be tailored gradually.

Finally, we give our sub-algorithm based on nonlinear center point to generate a good linear approximation for $X_{NL}$ in full details.

---

**Sub-algorithm: linear approximation (LA)**
0) Initialization: $\varepsilon^{LP} \geq 0$, $\Omega^R = \emptyset$, $k_{max}^{LP}$, $k^{LP} = 1$;
1) Center point step: solve (24) to obtain $\chi_{cp}^{NLP}$.
2) LP step: solve (25) to obtain $\chi^{LP}$.
3) If $g(\chi^{LP}) > \varepsilon^{LP}$, then conduct a line search between $\chi^{LP}$ and $\chi_{cp}^{NLP}$ to obtain $\hat{\chi}$, and go to step 4); otherwise, stop.
4) Let $\Omega^R = \Omega^R \cup \{\hat{\chi}\}$, update $X_{NL}^R(\Omega^R)$. If $k^{LP} \leq k_{max}^{LP}$, set $k^{LP} = k^{LP} + 1$ and go to step 2); else stop.

---

Note that, in this subsection, we don't take integer requirements into consideration. However, we obtain a *high-quality* overestimated polyhedral set for the feasible set of the relaxation of UC-CET. And along with the increasing of iteration times, $\chi^{LP}$ will gradually approach to optimal solution of the relaxation of UC-CET.

*2) Find integer solutions*

• Integer ellipsoid center point step

For giving points set $\Omega^R = \{\hat{\chi}^1, \hat{\chi}^2, \dots, \hat{\chi}^\hbar\}$ and $\Omega^F = \{\chi_{icp}^{u,1}, \chi_{icp}^{u,2}, \dots, \chi_{icp}^{u,m}\} \subseteq X_{NL} \cap X_L \cap X_I$. Let
$$X_{NL}^{R,B}(\Omega^R) = \left\{ (\chi, r) \left| \begin{array}{c} \tilde{g}_{\hat{\chi}^s}(\chi) + \mu r \|\nabla \tilde{g}_{\hat{\chi}^s}(\chi)\|_2^2 \leq 0 \\ \tilde{g}_{\hat{\chi}^s}(\chi) = \sum_{i=1}^N \sum_{t=1}^T [2c_i \hat{\tilde{P}}_{i,t}^s \tilde{P}_{i,t} - c_i u_{i,t} (\hat{\tilde{P}}_{i,t}^s)^2] - \eta_{CET} \\ s = 1, \dots, \hbar \end{array} \right. \right\} \tag{27}$$

$$X_{Obj}^{R,B}(\Omega^F) = \left\{ (\chi, r) \left| \ell(\chi) + \frac{1}{\mu} r \|\nabla \ell(\chi)\|_2^2 \leq \min_{i=1,\dots,m} \{\ell(\chi_{icp}^{u,i})\} \right. \right\}, \tag{28}$$

where $X_{NL}^{R,B}(\Omega^R)$ differ from $X_{NL}^{PR}(\Omega)$ in additional a new term $\mu r \|\nabla \tilde{g}_{\hat{\chi}^s}(\chi)\|_2^2$, and $\mu$ is a given ellipsoid factor.

**Definition 1** Let $X^{R,B} = X_{NL}^{R,B}(\Omega^R) \cap X_{Obj}^{R,B}(\Omega^F)$, *integer ellipsoid center* of $X^{R,B}$ can be found by solving the MILP
$$(\chi_{icp}, \hat{r}) = \arg \max_{\chi \in X_L \cap X_I \cap X^{R,B}, r \geq 0} r, \tag{29}$$
where $r$ is a variable and it can be viewed as the radius of the inscribed ellipsoid.

In this step, we want to find a center point (integer ellipsoid center) $\chi_{icp}$ of $X^{R,B}$ as a new trial solution by solving (29). However, $\chi_{icp}$ is not necessarily a feasible solution. There are two possibilities for $\chi_{icp}$.

If the trial solution $\chi_{icp}$ is feasible, then we need to decide the next step according to the value of $\hat{r}$. If $\hat{r} = 0$, then the CP algorithm stops and $\chi_{icp}$ is the optimal solution; Otherwise, we perform fixed-integer neighborhood search step.

If the trial solution $\chi_{icp}$ violates the nonlinear constraint, then we can perform integer variable feasibility adjustment step.

• Fixed-integer neighborhood search step

Now, we know that the integer variables $u_{icp}$ are feasible, but the trial solution $\chi_{icp}$ may still not be the best solution with the integer variables $u_{icp}$. we will therefore fix the integer variables in (21) to the values given by $u_{icp}$, resulting in the following convex quadratically constrained programming (QCP) problem
$$\chi_{icp}^u = \arg \min_{\chi \in X_{NL} \cap X_L \cap X_I, u = u_{icp}} \ell(\chi). \tag{30}$$
Then let $\Omega^F = \Omega^F \cup \{\chi_{icp}^u\}$, If $g(\chi_{icp}^u) = 0$, we need to let $\hat{\chi} = \chi_{icp}^u$ and update $\Omega^R = \Omega^R \cup \{\hat{\chi}\}$.

• Integer variable feasibility adjustment step

For the infeasible solution $\chi_{icp}$, we first need to verify the feasibility of the integer variables $u_{icp}$. The QCP subproblem that fixed the integer decision variables $u_{icp}$ is defined as
$$\min_{\chi \in X_L \cap X_I, u = u_{icp}} h$$
$$s.t. \, g(\chi) \leq h \tag{31}$$
$$h \geq 0,$$
where $h$ is a variable.

Obviously, if $h = 0$, then $u_{icp}$ is feasible and next perform fixed-integer neighborhood search step; if $h > 0$, then $u_{icp}$ is not feasible and next perform a line search between $\chi_{icp}$ and $\chi_{cp}^{NLP}$ to obtain $\hat{\chi}$, and update $\Omega^R = \Omega^R \cup \{\hat{\chi}\}$.

TABLE I
NUMBER OF UNITS IN EACH PROBLEM INSTANCE

| No. | Unit | | | | | | | | Total Units | No. | Unit | | | | | | | | Total Units |
|---|---|---|---|---|---|---|---|---|---|---|---|---|---|---|---|---|---|---|---|
| | 1 | 2 | 3 | 4 | 5 | 6 | 7 | 8 | | | 1 | 2 | 3 | 4 | 5 | 6 | 7 | 8 | |
| 1 | 12 | 11 | 0 | 0 | 1 | 4 | 0 | 0 | 28 | 12 | 58 | 50 | 15 | 7 | 16 | 18 | 7 | 12 | 183 |
| 2 | 13 | 15 | 2 | 0 | 4 | 0 | 0 | 1 | 35 | 13 | 55 | 48 | 18 | 5 | 18 | 17 | 15 | 11 | 187 |
| 3 | 15 | 11 | 0 | 1 | 4 | 5 | 6 | 3 | 45 | 14 | 240 | 220 | 0 | 0 | 20 | 80 | 0 | 0 | 560 |
| 4 | 10 | 10 | 2 | 5 | 7 | 5 | 6 | 5 | 50 | 15 | 260 | 300 | 40 | 0 | 80 | 0 | 0 | 20 | 700 |
| 5 | 13 | 12 | 5 | 7 | 2 | 5 | 4 | 6 | 54 | 16 | 300 | 260 | 40 | 120 | 60 | 20 | 20 | 60 | 880 |
| 6 | 46 | 45 | 8 | 0 | 5 | 0 | 12 | 16 | 132 | 17 | 300 | 220 | 0 | 20 | 80 | 100 | 120 | 60 | 900 |
| 7 | 40 | 54 | 14 | 8 | 3 | 15 | 9 | 13 | 156 | 18 | 300 | 260 | 60 | 140 | 100 | 60 | 40 | 20 | 980 |
| 8 | 51 | 58 | 17 | 19 | 16 | 1 | 2 | 1 | 165 | 19 | 200 | 200 | 40 | 100 | 140 | 100 | 120 | 100 | 1000 |
| 9 | 43 | 46 | 17 | 15 | 13 | 15 | 6 | 12 | 167 | 20 | 340 | 320 | 20 | 60 | 20 | 140 | 40 | 80 | 1020 |
| 10 | 50 | 59 | 8 | 15 | 1 | 18 | 4 | 17 | 172 | 21 | 240 | 340 | 80 | 140 | 100 | 40 | 0 | 100 | 1040 |
| 11 | 53 | 50 | 17 | 15 | 16 | 5 | 14 | 12 | 182 | 22 | 260 | 240 | 100 | 140 | 40 | 100 | 80 | 120 | 1080 |

**Algorithm: CP**
0) Initialization: $\varepsilon^r \geq 0$, $k^{\text{MILP}} = 0$, $\varepsilon^g \geq 0$, $\varepsilon^h \geq 0$, $\Omega^F = \emptyset$; Call Sub-algorithm LA to obtain $\Omega^R$ and $\chi_{\text{cp}}^{\text{NLP}}$.
1) Compute $\mu = \frac{1}{1+10^3 e^{-3k^{\text{MILP}}}}$, set $k^{\text{MILP}} := k^{\text{MILP}} + 1$.
2) Integer ellipsoid center point step: construct $X_{\text{NL}}^{\text{R,B}}(\Omega^R)$ according to (27), construct $X_{\text{Obj}}^{\text{R,B}}(\Omega^F)$ according to (28), and solve problem (29) to obtain $(\chi_{\text{icp}}, \hat{r})$.
3) If $g(\chi_{\text{icp}}) < \varepsilon^g$, go to step 4); otherwise, go to step 5).
4) If $\hat{r} < \varepsilon^r$, stop; otherwise, go to step 6).
5) Integer variable feasibility adjustment step: solve problem (31) to obtain $h$. If $h \geq \varepsilon^h$, then conduct a line search between $\chi_{\text{icp}}$ and $\chi_{\text{cp}}^{\text{NLP}}$ to obtain $\hat{\chi}$, update $\Omega^R = \Omega^R \cup \{\hat{\chi}\}$, and go to step 1); otherwise, go to step 6).
6) Fixed-integer neighborhood search step: solve problem (30) to obtain $\chi_{\text{icp}}^u$, and update $\Omega^F = \Omega^F \cup \{\chi_{\text{icp}}^u\}$. If $-\varepsilon^g \leq g(\chi_{\text{icp}}^u) \leq \varepsilon^g$, let $\hat{\chi} = \chi_{\text{icp}}^u$, update $\Omega^R = \Omega^R \cup \{\hat{\chi}\}$, and go to step 1); otherwise, directly go to step 1).

In this subsection, the constraint (28) plays an important role. On the one hand, it can reduce the search space and exclude all solutions that are worse than the current one. That is, once we get a feasible solution, it can force the algorithm to search for solutions with better objective function value than the feasible solution obtained. On the other hand, it will force $r$ to converge to zero. Because $r = 0$ means that $X^{R,B}$ has an empty interior. In other words, the algorithm can converge to the optimal solution.

## IV. NUMERICAL RESULTS AND ANALYSIS

To assess the performance of the proposed algorithm, we compared it with the state-of-the-art solver CPLEX. Twenty-two test instances with 28 to 1080 units running for a time span of 24h were used in our experiments. As presented in [21], these instances were created by replicating an eight-unit data set. Note that we removed the instances with repeated number of units and some small-scale instances. Table 1 shows the number of units in each instance. The parameters in CP algorithm are set as follows: $\varepsilon^{\text{LP}} = 0.001$, $k_{\max}^{\text{LP}} = 1000$, $\varepsilon^r = 0.001$, $\varepsilon^g = 0.001$, $L = 4$. All the codes and instances of the simulations for this paper can be freely downloaded from [34]. The machine on which we perform all of our computations is Dell XPS8930 with 16 GB of RAM and Intel i7-8700 3.2GHz CPU, running MS-Windows 10 and Matlab2016b. CPLEX 12.7.1 was used for solving LP, MILP, QCP and MIQCP problems.

The goals of this section are twofold. The first goal is to demonstrate the tightness and compactness of the linear approximate model generated by LA sub-algorithm. The second is that we will provide evidence to prove that: 1) CP algorithm can find high-quality feasible solutions faster than CPLEX and it is suitable to solve large-scale power systems; 2) Introducing integer ellipsoid center can improve the performance of CP algorithm.

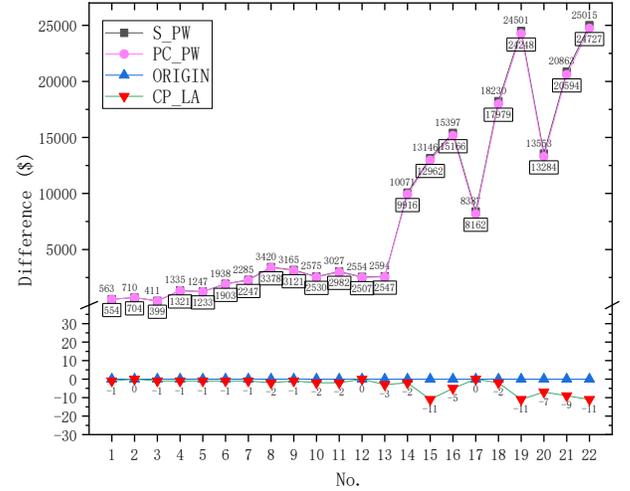

Fig.1. Comparison of the tightness of the four MIP formulations

### A. Simulation Results of LA sub-algorithm

For purpose of comparison, we build four different formulations, which are the linear approximate formulation generated by LA sub-algorithm, the standard piecewise linear approximate formulation, the piecewise linear approximate formulation with perspective-cut and the original formulation, respectively. Comparing the tightness of the four different MIP formulations can be achieved by comparing the optimal values of their continuous relaxation [18]. The larger the optimal value, the tighter the model. Therefore, the tightness of the models can be measured by the difference defined as $Z_{\text{CR\_ORIG}} - Z_{\text{CR\_MILP}}$. In this expression, $Z_{\text{CR\_MILP}}$ represents the optimal objective function value for the continuous relaxation of linear approximate MILP model for the original MINLP problem, and $Z_{\text{CR\_ORIG}}$ represents the optimal objective function value for the continuous relaxation of the original MINLP problem. To ensure fairness, the accuracy of CPLEX solver was set to the default value.

| TABLE II |
|---|
| THE OPTIMAL VALUES OF CONTINUOUS RELAXATION FOR THE ORIGINAL FORMULATION |

| No. | $Z_{CR\_ORIG}$ ($) | No. | $Z_{CR\_ORIG}$ ($) |
|---|---|---|---|
| 1 | 3729354 | 12 | 19647099 |
| 2 | 4833754 | 13 | 19258719 |
| 3 | 4678900 | 14 | 74531682 |
| 4 | 4319503 | 15 | 93778186 |
| 5 | 4959390 | 16 | 99962345 |
| 6 | 15432016 | 17 | 93580251 |
| 7 | 16799354 | 18 | 105329952 |
| 8 | 19661378 | 19 | 86369462 |
| 9 | 16961670 | 20 | 113801656 |
| 10 | 19034194 | 21 | 109887781 |
| 11 | 19219431 | 22 | 99194795 |

| TABLE III |
|---|
| THE NUMBER OF CUTTING PLANES GENERATED BY THREE LINEAR APPROXIMATE TECHNIQUES |

| No. | CP_LA | S_PW | PC_PW | No. | CP_LA | S_PW | PC_PW |
|---|---|---|---|---|---|---|---|
| 1 | 1 | 5 | 5 | 12 | 1 | 5 | 5 |
| 2 | 1 | 5 | 5 | 13 | 2 | 5 | 5 |
| 3 | 1 | 5 | 5 | 14 | 1 | 5 | 5 |
| 4 | 1 | 5 | 5 | 15 | 1 | 5 | 5 |
| 5 | 1 | 5 | 5 | 16 | 1 | 5 | 5 |
| 6 | 1 | 5 | 5 | 17 | 1 | 5 | 5 |
| 7 | 1 | 5 | 5 | 18 | 1 | 5 | 5 |
| 8 | 1 | 5 | 5 | 19 | 1 | 5 | 5 |
| 9 | 1 | 5 | 5 | 20 | 1 | 5 | 5 |
| 10 | 1 | 5 | 5 | 21 | 1 | 5 | 5 |
| 11 | 1 | 5 | 5 | 22 | 1 | 5 | 5 |

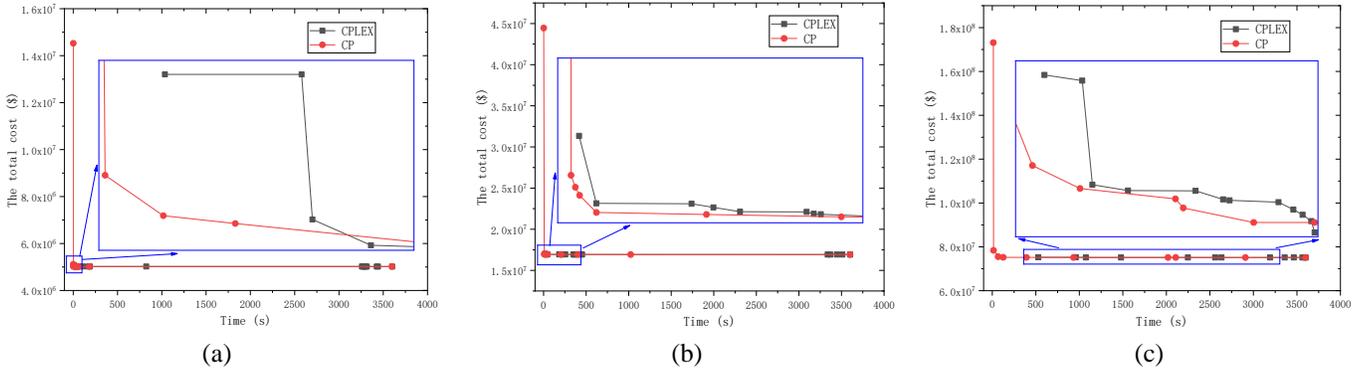

Fig.2. The comparison of the objective values of all feasible solutions found by CP algorithm and CPLEX within 3600s ((a)(b)(c) represents the instances of 54-,156- and 560-unit systems respectively).

In Fig.1, we show the difference values for all formulations for all test systems, where the difference for the original formulation itself is always equal to 0. The detailed values of difference for the three linear approximate formulations have been labeled in Fig.1 as well. In this figure, "No." denotes the No. of each problem instance. "S_PW" represents the continuous relaxation of the standard piecewise linear approximate formulation. "PC_PW" represents the continuous relaxation of the piecewise linear approximate formulation with perspective-cut. "ORIGIN" represents the continuous relaxation of the original formulation. "CP_LA" represents the continuous relaxation of the linear approximate formulation generated by LA sub-algorithm. The detailed optimal values of continuous relaxation for the original formulation are also given in Table II.

As shown in Fig.1, it can be seen from the instances of No. 1 – 14 that the tightness of ORIGIN and CP_LA models are similar. But From the remaining instances, except for No. 17, we can see the CP_LA model is tighter than ORIGIN model. And PC_PW model is tighter than S_PW model for all instances. Because the perspective-cut can remove the continuous-relaxed feasible area without affecting the feasible solutions [33]. Meanwhile, we can see in Fig.1 that the CP_LA model is tighter than two piecewise linear approximations. This is because the linearization points of piecewise linear approximation are chosen at regular intervals, and these points are likely to be far from the optimal solution, resulting in an unsatisfactory linear approximation. Consequently, the linear approximation generated by LA sub-algorithm is tighter than other three models.

Now, the compactness of the three linear approximate formulations needs to be compared. Table III presents the number of cutting planes generated by these linear approximate techniques. It can be seen from the table that the number of cut planes generated by LA sub-algorithm is much less than that generated by the piecewise linear approximate technique. This is because the linearization points generated by the LA sub-algorithm are closer to the optimal solution. Therefore, the linear approximate model generated by LA sub-algorithm is more compact than the other two piecewise approximate models.

### B. Simulation Results of CP Algorithm

To compare the performance of CP algorithm and CPLEX in finding feasible solutions, we compared the objective values of all feasible solutions found by CP algorithm and CPLEX within 3600s. In Fig. 2, (a)(b)(c) represents the comparison of the results for the instances of 54-,156- and 560-unit systems, respectively. we can see from the blue rectangular area that most of the red polylines are located below the black polyline, which proves that CP algorithm can find high-quality feasible solutions faster than CPLEX.

TABLE IV
THE RESULTS OBTAINED BY THE CP ALGORITHM AND CPLEX

| No. | $N$ | A solution within 5% of $Z_{CR\_ORIG}$ | | | | | A solution within 1% of $Z_{CR\_ORIG}$ | | | | |
|---|---|---|---|---|---|---|---|---|---|---|---|
| | | CP | | | CPLEX | | CP | | | CPLEX | |
| | | $F_{TC}$ ($) | $C_{time}$(s) | $N_{it}$ | $F_{TC}$ ($) | $C_{time}$ (s) | $F_{TC}$ ($) | $C_{time}$ (s) | $N_{it}$ | $F_{TC}$ ($) | $C_{time}$ (s) |
| 1 | 28 | 3872683 | **0.5** | 2 | 3774443 | 1.7 | 3766638 | **7.3** | 8 | 3766295 | 24.6 |
| 2 | 35 | 4877182 | 10 | 6 | 4889210 | **2.6** | 4877182 | 10 | 6 | 4881746 | **4** |
| 3 | 45 | 4747745 | 4.5 | 2 | 4740128 | **2.5** | 4723491 | 10 | 4 | 4724710 | **4.6** |
| 4 | 50 | 4391144 | 7.6 | 2 | 4374581 | **4.9** | 4368742* | 244* | 11 | 4368169* | 3281* |
| 5 | 54 | 5039960 | **1.4** | 2 | 5061476 | 4 | 5020375* | 182* | 10 | 5018696* | 3441* |
| 6 | 132 | 15729899 | **6** | 3 | 15430112 | 15 | 15566681 | **25** | 6 | 15579551 | 39 |
| 7 | 156 | 16977451 | **7.3** | 2 | 17020654 | 19 | 16964266 | **13** | 3 | 16946710 | 44 |
| 8 | 165 | 19861831 | **12** | 3 | 19832201 | 57 | 19855840 | **24** | 4 | 19832201 | 57 |
| 9 | 167 | 17150541 | **11** | 2 | 17176588 | 58 | 17129770 | **555** | 6 | 17130585 | 559 |
| 10 | 172 | 19316617 | **9** | 3 | 19319873 | 18 | 19208198 | **18** | 5 | 19223246 | 49 |
| 11 | 182 | 19433285 | **6** | 2 | 19430856 | 48 | 19410853 | **42** | 4 | 19401402 | 188 |
| 12 | 183 | 19837419 | **10** | 2 | 19824161 | 163 | 19837419 | **10** | 2 | 19824161 | 163 |
| 13 | 187 | 19483403 | **18** | 2 | 19472800 | 103 | 19439896 | **64** | 4 | 19439121 | 128 |
| 14 | 560 | 75520665 | **67** | 2 | 75241177 | 528 | 75217963 | **125** | 3 | 75241177 | 528 |
| 15 | 700 | 94989087 | **108** | 2 | 95838773 | 232 | 94657270 | **204** | 3 | 94528680 | 2286 |
| 16 | 880 | 101044912 | **150** | 3 | 100959878 | 1835 | 100965135* | 3361* | 5 | 100959878 | 1835 |
| 17 | 900 | 94582155 | **149** | 2 | 94466838 | 1697 | 94475429 | **301** | 3 | 94466838 | 1697 |
| 18 | 980 | 106558454 | **80** | 2 | 106359931 | 1159 | 106442881* | 1031* | 3 | 106359931 | 1159 |
| 19 | 1000 | 87668808 | **225** | 2 | 87359437 | 1044 | 87427528* | 2427* | 4 | 87350356* | 3565* |
| 20 | 1020 | 115068390 | **159** | 3 | 115301960 | 407 | 114927633 | **338** | 4 | 115051936* | 2488* |
| 21 | 1040 | 111143237 | **96** | 2 | 111256240 | 1297 | 111036268* | 2095* | 3 | 111256240* | 1297* |
| 22 | 1080 | 102185114 | **127** | 2 | - | - | 100602117* | 130* | 2 | - | - |

To further illustrate, we set two reference points, which are a solution within 5% of $Z_{CR\_ORIG}$ ($1.05 * Z_{CR\_ORIG}$) and a solution within 1% of $Z_{CR\_ORIG}$ ($1.01 * Z_{CR\_ORIG}$). The results of the CP algorithm and CPLEX are presented in Table IV. The column "$N$" denotes the total number of generating units. "$C_{time}$" and "$N_{it}$" denote the execution time and the number of iterations required to find a solution within 5% or 1% of $Z_{CR\_ORIG}$, respectively. The sign "*" indicates that the incumbent solution is the best one that the CP algorithm or CPLEX solver can find but still does not meet the requirements. The sign "-" indicates that no solution found within limited time. For CP algorithm, we solve the LP problem to 0.5% optimality and the NLP problem to 0.1% optimality using CPLEX. The execution time of CP algorithm and CPLEX is limited to 3600s.

In Table IV, for all instances, the CP algorithm is able to find a feasible solution within 5% of $Z_{CR\_ORIG}$ in less than 225 seconds. For some instances, a feasible solution within 1% of $Z_{CR\_ORIG}$ can be obtained by CP algorithm in less than 560 seconds. Compared to CPLEX, the proposed algorithm has obvious advantages. At the same time, it can be seen from Table IV that some results did not meet the requirements (these results were marked by *). There are two main reasons causing this situation. On the one hand, it is difficult for CP algorithm and CPLEX to find a solution within 1% of $Z_{CR\_ORIG}$ owing to the value of $Z_{CR\_ORIG}$ is too small, such as 50-, 54-, 1000-, and 1040-unit systems. On the other hand, different instances are suitable for different algorithms, for example, CPLEX can find a better solution than CP algorithm for 880- and 980-unit systems, but for 1020-unit systems, the opposite is true. In general speaking, CP algorithm has relatively satisfactory performance in comparison with CPLEX.

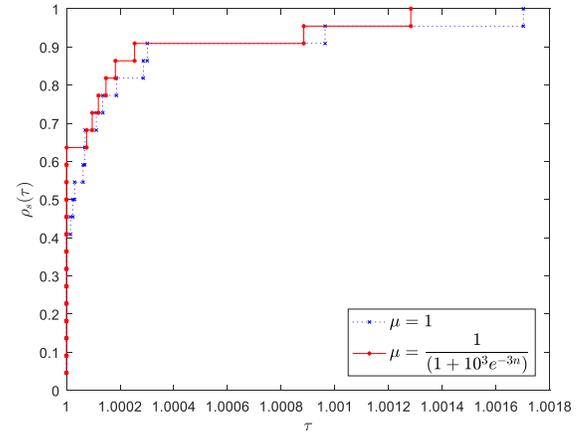

Fig.3. Performance profile on the objective value for CP algorithm when using $\mu = \frac{1}{1+10^3 e^{-3n}}$ ($n = k^{MILP}$) and $\mu = 1$

Meanwhile, in order to present the performance of CP algorithm that introduces the integer ellipsoid center, we adopt the performance profiles [35] to compare the performance according to the objective value of the best solution obtained by the CP algorithm when using different values of $\mu$. Then we define $t_{p,s}$ as performance of method $s$ while solving problem $p$,

$$\rho_s(\tau) = \frac{1}{|\mathcal{P}|} size\{p \in \mathcal{P} : \frac{t_{p,s}}{min\{t_{p,s}:s \in \mathcal{S}\}} \leq \tau\} \quad (32)$$

where $\mathcal{S}$ is the set of methods and $\mathcal{P}$ is the set of problems [35].

The comparison of the results is shown in Fig.3. As can be seen from the figure, the solid red line corresponds to $\mu = \frac{1}{1+10^3 e^{-3n}} (n = k^{\text{MILP}})$ is mostly located above the blue dotted line corresponds to $\mu = 1$, which means that introduction of the integer ellipsoid center in CP algorithm will improve its performance. In other word, the CP algorithm introducing the integer ellipsoid center can find better solutions.

## V. CONCLUSION

In this paper, a deterministic global optimization method has been presented to solve the UC-CET. By constructing perspective-cuts at linearization points acquired with solving a sequence of linear continuous-relaxed subproblems and performing line search procedure, the proposed method can quickly generate a tight linear approximation of UC-CET. And then the method can rapidly obtain high-quality solutions by iteratively finding the integer ellipse center of the current linear approximation, searching for the neighborhood of these points and adding new linear constraints. We compared the performance of the CP algorithm with the performance of CPLEX. Numerical results have demonstrated the CP algorithm can find high-quality solutions faster.